\documentclass[12pt]{amsart}

\usepackage{amsmath}
\usepackage{amssymb}
\usepackage{pdfsync}

\newcommand{\C}{{\mathbb C}}       
\newcommand{\R}{{\mathbb R}}       
\newcommand{\Z}{{\mathbb Z}}       

\newcommand{\FF}{{\mathcal F}}
\newcommand{\HH}{{\mathcal H}}
\newcommand{\LL}{{\mathcal L}}

\newcommand{\BZ}{{\mathcal B}}
\newcommand{\GZ}{{\mathcal G}}

\newcommand{\CC}{{\mathcal C}}

\newcommand{\diam}{{\rm diam}}
\newcommand{\dist}{{\rm dist}}

\newcommand{\fiproof}{{\hspace*{\fill} $\square$ \vspace{2pt}}}

\newcommand{\rf}[1]{{(\ref{#1})}}

\newcommand{\supp}{{\rm supp}}

\newcommand{\vphi}{{\varphi}}
\newcommand{\ve}{{\varepsilon}}
\newcommand{\vv}{{\vspace{2mm}}}

\newcommand{\wt}[1]{{\widetilde{#1}}}
\newcommand{\wh}[1]{{\widehat{#1}}}

\newcommand{\noi}{\noindent}

\newtheorem{theorem}{Theorem}[section]
\newtheorem{lemma}[theorem]{Lemma}

\newtheorem*{theorem*}{Theorem}

\theoremstyle{definition}

\newtheorem{example}[theorem]{Example}

\theoremstyle{remark}
\newtheorem{rem}[theorem]{Remark}

\numberwithin{equation}{section}

\newcommand{\brem}{\begin{rem}}
\newcommand{\erem}{\end{rem}}

\newcommand{\bexam}{\begin{example}}
\newcommand{\eexam}{\end{example}}
\newcommand{\sss}{{\rm Stop}}

\begin{document}

\title[Calder\'on-Zygmund capacities and Wolff potentials]{
Calder\'on-Zygmund capacities and Wolff potentials on Cantor sets}

\author[X. Tolsa]{Xavier Tolsa}
\address{Instituci\'{o} Catalana de Recerca i Estudis Avan\c{c}ats (ICREA) and Departament de Matem\`{a}tiques,  Universitat Aut\`{o}noma de Barcelona, 08193 Bellaterra (Barcelona), Catalonia}
\email{{\tt xtolsa@mat.uab.cat}}
\urladdr{http://mat.uab.cat/~xtolsa}

\thanks{Partially supported by grant MTM2007-62817 (Spain) and
2009-SGR-420 (Catalonia). Part of this paper was done during the attendance to the research semester
``Harmonic Analysis, Geometric Measure Theory and Quasiconformal Mappings'' at the CRM (Barcelona), in 2009.}

\begin{abstract}
We show that, for some Cantor sets in $\R^d$, 
the capacity $\gamma_s$ associated to the 
$s$-dimensional Riesz kernel $x/|x|^{s+1}$ is comparable to the capacity
$\dot C_{\frac23(d-s),\frac32}$ from non linear potential theory. 
It is an open problem to show that, when $s$ is positive and non integer, they are comparable for all compact sets in $\R^d$. We also discuss other open questions in the area.
\end{abstract}

\maketitle

\section{Introduction}

In the first part of this paper we show that, for some Cantor sets in $\R^d$, 
the capacity $\gamma_s$ associated to the 
$s$-dimensional Riesz kernel $x/|x|^{s+1}$ is comparable to the capacity
$\dot C_{\frac23(d-s),\frac32}$ from non linear potential theory. 
It is an open problem to show that, when $s$ is a positive and non integer, they are comparable for all compact sets in $\R^d$. In the last part of the paper, we discuss other related open questions.

To state our results in detail we need to introduce some notation. For $0<s<d$, 
the $s-$dimensional Riesz kernel is defined by
$$K^s(x)=\frac{x}{|x|^{s+1}},\quad x\in \R^d,  \, x\neq0.$$
Notice that this is a vectorial kernel.  The
$s-$dimensional Riesz transform (or $s-$Riesz transform) of a real
Radon measure $\nu$ with compact support is
 $$R^s \nu(x) =\int
K^s(y-x)\, d\nu(y),\qquad x\not\in\supp(\nu).$$ 
Although the preceding integral converges a.e.\ with respect to Lebesgue measure, the convergence may fail for
$x\in\supp(\nu)$. This is the reason why one considers the truncated
$s-$Riesz transform of $\nu$, which is defined as $$R^s_\ve \nu(x)
=\int_{|y-x|>\ve} K^s(y-x)\, d\nu(y),\qquad x\in\R^d,\,\ve>0.$$
These definitions also make sense if one consider distributions
instead of measures. Given a compactly supported distribution $T$, set
$$R^s(T) = K^s*T$$
 (in the principal value sense for $s=d$), and analogously $$R^s_\ve(T) =
K^s_\ve*T,$$ where $K^s_\ve(x) = \chi_{|x|>\ve}\,x/|x|^{s+1}$.

Given a positive Radon measure with compact support and a function
$f\in L^1(\mu)$, we consider the operators $R^s_\mu(f) :=
R^s(f\,d\mu)$ and $R^s_{\mu,\ve}(f) := R^s_\ve(f\,d\mu)$. We say
that $R^s_\mu$ is bounded on $L^2(\mu)$ if $R^s_{\mu,\ve}$ is
bounded on $L^2(\mu)$ uniformly in $\ve>0$, and we set
$$\|R^s_\mu\|_{L^2(\mu)\to L^2(\mu)} =
\sup_{\ve>0}\|R^s_{\mu,\ve}\|_{L^2(\mu)\to L^2(\mu)}.$$ 

Given a compact set $E\subset \R^d$, the capacity $\gamma_s$
of $E$ is
\begin{equation} \label{defkap}
\gamma_s(E) = \sup|\langle T,1\rangle|,
\end{equation}
where the supremum is taken over all distributions $T$ supported
on $E$ such that $\|R^{s}(T)\|_{L^\infty(\R^d)} \leq 1$. 
Following \cite{Volberg}, we call $\gamma_s$ the $s$-dimensional Calder\'on-Zygmund capacity. The case $s=d-1$ is particularly relevant: 
$\gamma_{d-1}$ coincides with the
capacity $\kappa$ introduced by Paramonov \cite{Paramonov}
 in order to
study problems of $\CC^1$ approximation by harmonic functions in
$\R^d$ (the reader should notice that $\kappa$ is called $\kappa'$
in \cite{Paramonov}). When $d=2$ and $s=1$, $z/|z|^{s+1}$
coincides with the complex conjugate of the Cauchy kernel $1/z$. Thus, if
one allows $T$ to be a complex distribution in the supremum above, then 
$\gamma_1$ is the analytic capacity.

If we restrict the supremum in \rf{defkap} to distributions
$T$ given by positive Radon measures supported on $E$, we obtain
the capacities $\gamma_{s,+}$. Clearly, we have $\gamma_s(E)\geq \gamma_{s,+}(E)$.
On the other hand, the opposite inequality also holds (up to a multiplicative absolute constant $c_s$):
$$\gamma_s(E)\leq c_s\,\gamma_{s,+}(E).$$
This was first shown for $s=1,d=2$ by the author \cite{Tolsa-sem}, and it was extended to
the case $s=d-1$ by Volberg \cite{Volberg}. For other values of $s$, this can be proved
by combining the techniques from \cite{Volberg} with others from 
\cite{Mateu-Prat-Verdera} (see \cite{Prat-pers}).

Now we turn to non linear potential theory. Given
$\alpha>0$ and $1<p<\infty$ with $0<\alpha p <2$, the capacity $\dot C_{\alpha,p}$
of $E\subset\R^d$ is defined as
$$\dot C_{\alpha,p}(E) = \sup_\mu \mu(E)^p,$$
where the supremum runs over all positive measures $\mu$ supported on $E$ such that
$$I_\alpha(\mu)(x) = \int \frac1{|x-y|^{2-\alpha}}\,d\mu(x)$$
satisfies $\|I_\alpha(\mu)\|_{p'}\leq 1$, where as usual $p'=p/(p-1)$. 

For our purposes, the characterization of $\dot C_{\alpha,p}$ in terms of Wolff
potentials is more useful than its definition above. Consider
$$\dot W^\mu_{\alpha,p}(x) = \int_0^\infty \biggl(\frac{\mu(B(x,r))}{r^{2-\alpha p}}\biggr)^{p'-1}\,\frac{dr}r.$$
A well known theorem of Wolff asserts that
\begin{equation}\label{eqwoo}
\dot C_{\alpha,p}(E) \approx \sup_\mu \mu(E),
\end{equation}
where the supremum is taken over all measures $\mu$ supported on $E$ such that
$\dot W_{\alpha,p}^\mu(x)\leq 1$ for all $x\in E$ (see \cite[Chapter 4]{adamshedberg},
for instance).
The notation $A\approx B$ means that there is an absolute constant $c>0$, or depending
on $d$ and $s$ at most,
such that $c^{-1}A \leq B\leq cB$.

Mateu, Prat and Verdera showed in \cite{Mateu-Prat-Verdera} that if $0<s<1$, then
$$\gamma_s(E) \approx \dot C_{\frac23(d-s),\frac32}(E).$$ Notice that the
Wolff's potential for the capacity $ \dot C_{\frac23(d-s),\frac32}$ is 
$$\dot W^\mu_{\frac23(d-s),\frac32}(x) = \int_0^\infty \biggl(\frac{\mu(B(x,r))}{r^s}\biggr)^2\,\frac{dr}r.$$
When $s=1$ and $d=2$, from the characterization
of $\gamma_{1,+}$ in terms of curvature of measures, one easily 
gets $\gamma_{1}(E) \gtrsim \dot C_{\frac23,\frac32}(E)$.
Using analogous arguments (involving a symmetrization of the kernel and the $T(1)$ theorem), in \cite{Ei-Na-Vo} it has been shown that this also holds for all indices 
$0<s<d$:
$$\gamma_s(E) \gtrsim \dot C_{\frac23(d-s),\frac32}(E),$$ 
for any compact set $E\subset\R^d$.
The opposite inequality is false when $s$ is integer (for instance, if $E$ is 
contained in an $s$-plane and has positive $s$-dimensional Hausdorff measure, then
$\gamma_s(E)>0$, but $\dot C_{\frac23(d-s),\frac32}(E)=0$).
When $0<s<d$ is non integer, it is an open problem to prove (or disprove) that
$$\gamma_s(E) \lesssim \dot C_{\frac23(d-s),\frac32}(E).$$
See Section \ref{secopen} for more details and related questions.

In the present paper we show that the comparability
$\gamma_s(E) \approx \dot C_{\frac23(d-s),\frac32}(E)$ holds for some Cantor sets $E\subset\R^d$, which
are defined as follows. Given a sequence
$\lambda=(\lambda_n)^\infty_{n=1}$, $0\leq\lambda_n < 1/2$, we
construct $E$ by the following algorithm. Consider the
unit cube $Q^0=[0,1]^d$. At the first step we take $2^d$ closed
cubes inside $Q^0$, of side length $\ell_1=\lambda_1$, with sides
parallel to the coordinate axes, such that each cube contains a
vertex of $Q^0$. At the second step $2$ we apply the preceding procedure to
each of the $2^d$ cubes produced at step 1, but now using the
proportion factor $\lambda_2$. Then we obtain $2^{2d}$ cubes of
side length $\ell_2=\lambda_1\lambda_2$. Proceeding inductively,
we have at the $n-$th step $2^{nd}$ cubes $Q^n_j$, $1\leq j\leq
2^{nd}$, of side length $\ell_n=\prod_{j=1}^n\lambda_j$. We
consider $$
E_n=E(\lambda_1,\dots,\lambda_n)=\bigcup_{j=1}^{2^{nd}}Q_j^n, $$
and we define the Cantor set associated to
$\lambda=(\lambda_n)^\infty_{n=1}$ as $$
E=E(\lambda)=\bigcap_{n=1}^\infty E_n. $$ For example, if
$\lim_{n\to\infty}\ell_n / 2^{-nd/s}=1$, then the Hausdorff
dimension of  $E(\lambda)$ is $s$. If moreover
$\ell_n=2^{-nd/s}$ for each $n$, then
$0<\HH^s(E(\lambda))<\infty$, where $\HH^s$ stands for the
$s-$dimensional Hausdorff measure.
In the planar case ($d=2$), This class of Cantor sets first appeared in 
\cite{Garnett-lnm} (as far as we know),
and its study has played a very important role in the last advances concerning 
analytic capacity.

Our result reads as follows.

\begin{theorem} \label{teokappa}
Assume that, for all $n$, $0<
\lambda_n\le\tau_0<\frac{1}{2}.$ Denote $\theta_n=2^{-nd}/{\ell_n}^{\!\!s}$.
For any $N=1,2,\ldots$ we have
\begin{equation*}
\gamma_s(E_N)\approx \dot C_{\frac23(d-s),\frac32}(E_N) \approx
\biggl(\sum_{n=1}^{N}
\theta_n^{\,2}\biggr)^{-1/2},
\label{propo1}
\end{equation*}
where the constants involved in the relationship $\approx$ depend on $d$, $s$ and $\tau_0$, but not on $N$.
\end{theorem}

Observe that if $\mu$ is for the probability measure on $E_N$
 given by $\mu= \frac{\LL^d|E_N}{\LL^d(E_N)},$ where $\LL^d$
stands for the Lebesgue measure in $\R^d$, then 
$\theta_n=\mu(Q^n_j)/{\ell_n}^{\!\!s}$. So $\theta_n$ is the $s$-dimensional density
of $\mu$ on a cube from the $n$-th generation.

Showing that $\dot C_{\frac23(d-s),\frac32}(E_N) \approx
\Bigl(\sum_{n=1}^{N}
\theta_n^{\,2}\Bigr)^{-1/2}$ is not difficult, using the characterization of 
$\dot C_{\frac23(d-s),\frac32}$ in terms of Wolff's potentials (see Section 
\ref{secwo}). The difficult part of the theorem consists in showing that
\begin{equation}\label{eqgamss}
\gamma_s(E_N) \approx
\Bigl(\sum_{n=1}^{N} \theta_n^{\,2}\Bigr)^{-1/2}.
\end{equation}
 The main step in proving this result consists
in estimating the $L^2(\mu)$ norm of the $s$-dimensional Riesz transform $R_\mu^s$.

Let us remark that \rf{eqgamss} has been proved for analytic capacity ($s=1$, $d=2$) 
in \cite{MTV1} (using previous results from Mattila \cite{Mattila-cantors} and
Eiderman \cite{Eiderman}).
The arguments in \cite{MTV1} (as well as the ones in \cite{Mattila-cantors} and
Eiderman \cite{Eiderman}) rely heavily on the relationship between 
the Cauchy transform and curvature of measures. See \cite{Melnikov-curvatura} and \cite{Melnikov-Verdera} 
for more details on this relationship. 

In the case $s=d-1$, the
comparability \rf{eqgamss} was proved by Mateu and the author \cite{Mateu-Tolsa} under the 
additional assumption that $\lambda_n\geq 2^{-d/s}$ for all $n$, which is 
equivalent to saying that the sequence $\{\theta_n\}$ is non increasing.
It is not difficult to show that the arguments in \cite{Mateu-Tolsa} extend to all indices $0<s<d$.
However, getting rid of the assumption $\lambda_n\geq 2^{-d/s}$ is much more
 delicate. This is what
we carry out in this paper.

Let us also mention that in \cite{Garnett-Prat-Tolsa} it was shown that
the estimate \rf{eqgamss} also holds if one replaces $E_N$ by some bilipschitz image of itself,
also under the assumption $\lambda_n\geq 2^{-d/s}$. 
On the other hand, recently in \cite{Ei-Na-Vo} some examples of random Cantor
sets where the comparability $\gamma_s\approx \dot C_{\frac23(d-s),\frac32}$ holds 
have been studied.

The plan of the paper is the following. In Section \ref{secwo} we show that
$\dot C_{\frac23(d-s),\frac32}(E_N) \approx
\Bigl(\sum_{n=1}^{N}
\theta_n^{\,2}\Bigr)^{-1/2}$. The proof of \rf{eqgamss} is contained
in Sections \ref{secprelim}, \ref{secl21}, and \ref{secl22}. 
In the final Section \ref{secopen} we discuss open problems in connection
with Calder\'on-Zygmund capacities, Riesz transforms, and Wolff potentials.

Throughout all the paper, the letters $c,C$ will stand for
absolute constants (which may depend on $d$ and $s$) that may
change at different occurrences. Constants with subscripts, such
as $C_1$, will retain their values, in general.


\section{Proof of 
$\dot C_{\frac23(d-s),\frac32}(E_N) \approx
\Bigl(\sum_{n=1}^{N}
\theta_n^{\,2}\Bigr)^{-1/2}$}\label{secwo}

The proof of this result is essentially contained in \cite[Section 5.3]{adamshedberg}. However,
for the reader's convenience we give a simple and almost self-contained proof.

Recall that $\mu$ stands for the probability measure on $E_N$ defined by
$\mu= \frac{\LL^d|E_N}{\LL^d(E_N)}$.
Given $x\in E_N$, let $Q^n(x)$ denote the cube $Q^n_j$ from
the $n$-th generation in the construction of $E_N$ that contains $x$, so that $\ell(Q^n(x))=\ell_n$ is its 
side length.
It is straightforward to check that for all $x\in E_N$,
$$\dot W^\mu_{\frac23(d-s),\frac32}(x) = \int_0^\infty \biggl(\frac{\mu(B(x,r))}{r^s}\biggr)^2\,\frac{dr}r
\approx \sum_{n\geq0}\biggl(\frac{\mu(Q^n(x))}{\ell(Q^n(x))^s}\biggr)^2 = 
\sum_{n\geq0}\theta_n^2.$$
Thus, if we consider the measure 
$$\nu=\biggl(\sum_{n\geq0}\theta_n^2\biggr)^{-1/2}\,\mu,$$
we have $\dot W^\nu_{\frac23(d-s),\frac32}(x) \lesssim 1$ for all $x\in E_N$.
From \rf{eqwoo} we infer that
$$\dot C_{\frac23(d-s),\frac32}(E_N) \gtrsim \nu(E_N) = 
\biggl(\sum_{n\geq0}\theta_n^2\biggr)^{-1/2}.$$

To prove the converse inequality, we recall that given any Borel measure $\sigma$ on $\R^d$, for
any capacity $\dot C_{\alpha,p}$,
$$\dot C_{\alpha,p}\bigl(\{x\in\R^d:\,W_{\alpha,p}^\sigma(x)>\lambda\}\bigr)\leq c_{\alpha,p}\,\frac{\sigma(\R^d)}
{\lambda^{p-1}},\quad\mbox{for all $\lambda>0$.}$$
See Proposition 6.3.12 of \cite{adamshedberg}. If we apply this estimate to $\dot C_{\frac23(d-s),\frac32}$,\,
$\sigma=\mu$, and $\lambda \approx \sum_{n\geq0}\theta_n^2$, we get
$$\dot C_{\frac23(d-s),\frac32}(E_N)\leq 
\dot C_{\frac23(d-s),\frac32}\bigl(\{x\in\R^d:\,W_{\frac23(d-s),\frac32}^\sigma(x)> \lambda\}\bigr)
\lesssim \frac1{\Bigl(\sum_{n\geq0}\theta_n^2\Bigr)^{1/2}}.$$


\section{Preliminaries for the proof of $\gamma_s(E_N) \approx
\Bigl(\sum_{n=1}^{N}
\theta_n^{\,2}\Bigr)^{-1/2}$} \label{secprelim}

To simplify notation, to denote the $s$-dimensional Riesz transform of $\mu$ we will write $R\mu$ instead of 
$R^s\mu$, and also $K(x)$ instead of $K^s(x)=x/|x|^{s+1}$. Moreover, $\| \, \cdot \, \|$ stands for
the $L^2(\mu)$ norm.

Arguing as in \cite[Lemma 4.2]{Mateu-Tolsa}, it turns out that
 the estimate 
\begin{equation}\label{eqgamss2}
\gamma_s(E_N) \approx
\Bigl(\sum_{n=1}^{N}
\theta_n^{\,2}\Bigr)^{-1/2}
\end{equation}
 follows from the next result.

\begin{theorem} \label{difi0}
Let $\mu$ be the preceding probability measure supported on $E_N$.  We have $$\|R\mu\|^2 \approx \sum_{j=0}^N\theta_j^2.$$
\end{theorem}

\medskip
We will skip the arguments that show that \rf{eqgamss2} can be deduced from this theorem, which
the interested reader can
find in the aforementioned reference.
\medskip

Sections \ref{secl21} and \ref{secl22} of this paper are devoted to the proof of Theorem \ref{difi0}. In the remaining part of 
the current section, we introduce some additional notation that we will use below, and we prove a technical estimate.

Denote $\wt\Delta = \left\{Q^n_j: n\geq0,\,1\leq j\leq 2^{nd}\right\}$, where the $Q^n_j$'s are the cubes
which appear in the construction of the $E(\lambda)$.
Let $\Delta_n$ be the family of cubes in $\wt\Delta$ from the $n$-th generation. 
That is,
$\Delta_n= \{Q^n_j\}_{j=1}^{2^{nd}}$. For a fixed $N\geq1$, we set $\Delta =\bigcup_{n=1}^N\Delta_n$ 
(so $E_N$ is constructed using the cubes from $\Delta_N$).
 
 Given a cube $Q\subset\R^d$, we set
$$\theta(Q):=\frac{\mu(Q)}{\ell(Q)^s},$$
i.e.\ $\theta(Q)$ is the average
$s$-dimensional density of $\mu$ over $Q$. Thus
$\theta_n=\theta(Q)$ if $Q\in\Delta_n$.

 Given a cube $Q\in\Delta$ and a function $f\in L^1_{loc}(\mu)$, we define
$$S_Q f(x) = \frac1{\mu(Q)}\int_Q f\,d\mu \,\,\chi_Q(x).$$
Also, for $0\leq j\leq N$, we set
$S_j f = \sum_{Q\in\Delta_j} S_Q f.$
If we denote by $\FF(Q)$ the cubes from $\Delta$ which are sons of $Q$, we set
$$D_Q f(x) = \sum_{P\in\FF(Q)} S_P f(x) - S_Qf(x),$$
and for $0\leq j\leq N$ we denote $D_j f = \sum_{Q\in\Delta_j} D_Q f= S_{j+1}f-S_jf$.

Let $\Delta^0 = \Delta\setminus \Delta_N$. Notice that the functions $D_Qf$ and $D_Pf$ are orthogonal for 
$P\neq Q$. If $\int f\,d\mu=0$, then 
$$S_N f = \sum_{j=0}^{N-1} D_j f = \sum_{Q\in\Delta^0} D_Q f,$$
and thus
$$\|f\|^2\geq \|S_N f\|^2= \sum_{Q\in \Delta^0} \|D_Q f\|^2.$$
In particular, if we take $f= R\mu$, by antisymmetry $\int R\mu\,d\mu=0$, and thus
\begin{equation}\label{eqort1}
\|R\mu\|^2 \geq \|S_N (R\mu)\|^2 = \sum_{Q\in \Delta^0} \|D_Q (R\mu)\|^2.
\end{equation}

 Given cubes $Q,R\in\Delta$, 
we denote
\begin{equation}\label{defpqr}
p(Q) := \sum_{P\in\Delta:Q\subset P}\theta(P)\frac{\ell(Q)}{\ell(P)},\qquad
p(Q,R) := \sum_{P\in\Delta:Q\subset P\subset R}\theta(P)\frac{\ell(Q)}{\ell(P)}.
\end{equation}
For $0\leq j\leq N$, we denote $p_j := p(Q)$, for $Q\in\Delta_j$.

\begin{lemma}\label{lemnab}
Let $Q\in\Delta$ and $x,x'\in Q$. Let $\wh Q$ the parent of $Q$. Then we have
$$\bigl|R(\chi_{\R^d\setminus Q}\mu)(x) - R(\chi_{\R^d\setminus
Q}\mu)(x')\bigr| \leq C_1\,\frac{\ell(Q)}{\ell(\wh Q)}\,p(\wh Q).$$
Thus,
$$\bigl|R(\chi_{\R^d\setminus Q}\mu)(x) - R(\chi_{\R^d\setminus
Q}\mu)(x')\bigr| \leq C_1\,p(\wh Q)\leq C_2\,p(Q).$$
\end{lemma}

\begin{proof}
 We have
\begin{align*}
\bigl|R(\chi_{\R^d\setminus Q}\mu)(x) -  R(\chi_{\R^d\setminus Q}\mu)&(x')\bigr|\\
 & \leq \int_{\R^d\setminus Q}
|K(x-y)-K(x'-y)|\,d\mu(y)\\
& \leq C|x-x'|\int_{\R^d\setminus Q} \frac1{|x-y|^{s+1}}\,d\mu(y)\\
& \leq C|x-x'| \sum_{P\in\Delta:\,Q\subsetneq P}
\frac{\mu(P)}{\ell(P)^{s+1}} \leq C\, \frac{\ell(Q)}{\ell(\wh Q)}\,p(\wh Q).
\end{align*}
\end{proof}


\section{Proof of $\|R\mu\|^2 \lesssim \sum_{j=0}^N\theta_j^2.$}\label{secl21}

\begin{lemma}\label{lemdes11}
If $Q\in\Delta^0$ and $P$ is a son of $Q$, then
\begin{equation}\label{eqson1}
|S_P(R\mu) - S_Q(R\mu)|\lesssim p(Q).
\end{equation}
As a consequence,
\begin{equation}\label{eqson2}
\|D_Q (R\mu)\|^2\lesssim p(Q)^2\,\mu(Q).
\end{equation}
\end{lemma}

\begin{proof}
It is clear that \rf{eqson2} follows from \rf{eqson1}. To prove \rf{eqson1}, we use 
the antisymmetry of the kernel $K(x)$:
\begin{align}\label{eqdec0}
S_P(R\mu) - S_Q(R\mu)  &= 
S_P(R(\chi_{\R^d\setminus P}\mu)) - S_Q(R(\chi_{\R^d\setminus Q}\mu)) \nonumber\\
& = S_P(R(\chi_{Q\setminus P}\mu)) + 
S_P(R(\chi_{\R^d\setminus Q}\mu)) - S_Q(R(\chi_{\R^d\setminus Q}\mu)).
\end{align}
From Lemma \ref{lemnab} it follows that
$$|S_P(R(\chi_{\R^d\setminus Q}\mu)) - S_Q(R(\chi_{\R^d\setminus Q}\mu))|
\lesssim p(Q).$$
To estimate $S_P(R(\chi_{Q\setminus P}\mu))$ we take into account that
$\dist(Q\cap E_N\setminus P,P)\approx \ell(Q)$, and so for every  $x\in P$,
$$|R(\chi_{Q\setminus P}\mu)(x)|\lesssim\frac{\mu(Q)}{\ell(Q)^s} = \theta(Q)\leq p(Q).$$
From the preceding estimates and \rf{eqdec0}, we get \rf{eqson1}.
\end{proof}

\begin{lemma}\label{lemfa1}
We have
$$\|S_N (R\mu)\|^2 \lesssim \sum_{j=0}^{N-1}\theta_j^2\qquad\mbox{and}
\qquad\|R\mu\|^2 \lesssim \sum_{j=0}^N\theta_j^2.$$
\end{lemma}

\begin{proof}
By \rf{eqort1} and Lemma \rf{lemdes11},
$$\|S_N (R\mu)\|^2 = \sum_{Q\in \Delta^0} \|D_Q (R\mu)\|^2 \lesssim
\sum_{Q\in \Delta^0} p(Q)^2\mu(Q) = \sum_{j=0}^{N-1}p_j^2.$$
On the other hand, by Lemma \ref{lemnab}, for each $Q\in\Delta_N$ and $x\in Q$,
$$|S_Q (R\mu) - R(\chi_{\R^d\setminus Q}\mu)(x)| = 
|S_N (R\chi_{\R^d\setminus Q}\mu)) - R(\chi_{\R^d\setminus Q}\mu)(x)| \lesssim p(Q).$$
Using also that
$$\|\chi_Q R(\chi_Q \mu)\| \leq\theta(Q)\,\mu(Q)^{1/2},$$
we obtain
\begin{align*}
\|R\mu\|^2 & = \sum_{Q\in\Delta_N} \|\chi_Q R(\mu)\|^2 \leq 
2 \sum_{Q\in\Delta_N} \Bigl(\|\chi_Q R(\chi_Q\mu)\|^2 +
 \|\chi_Q R(\chi_{\R^d\setminus Q}\mu)\|^2\Bigr)
\\
&\leq2 \sum_{Q\in\Delta_N} \Bigl(\|\chi_Q R(\chi_Q \mu)\|^2 + \|R(\chi_{\R^d\setminus Q}\mu)
- S_N (R\chi_{\R^d\setminus Q}\mu))\|^2 + \|S_N (R\mu)\|^2\Bigr) \\
& \lesssim  \sum_{Q\in\Delta_N} \theta(Q)^2\mu(Q) + \sum_{Q\in\Delta_N} p(Q)^2\mu(Q)
+ \sum_{Q\in \Delta^0} p(Q)^2\mu(Q)\\
& \lesssim \sum_{Q\in \Delta} p(Q)^2\mu(Q) = \sum_{j=0}^N p_j^2.
\end{align*}
It only remains to show that
$\sum_{j=0}^M p_j^2\lesssim \sum_{j=0}^M \theta_j^2$ both for $M=N-1$ and $M=N$.
This follows easily from the definition of $p_j$ and Cauchy-Schwartz:
\begin{align} \label{eqpjtj}
\sum_{j=0}^M p_j^2 & = \sum_{j=0}^M 
\biggl( \sum_{k=0}^j \theta_k\,\frac{\ell_j}{\ell_k}\biggr)^2 \leq 
 \sum_{j=0}^M \biggl( \sum_{k=0}^j \theta_k^2\,\frac{\ell_j}{\ell_k}\biggr)
 \biggl( \sum_{k=0}^j \frac{\ell_j}{\ell_k}\biggr) \nonumber\\
 & \leq 2\sum_{j=0}^M \sum_{k=0}^j \theta_k^2\,\frac{\ell_j}{\ell_k}
 = 2 \sum_{k=0}^M \theta_k^2\sum_{j=k}^M \frac{\ell_j}{\ell_k} \leq
 4\sum_{k=0}^M \theta_k^2.
\end{align}
\end{proof}


\section{Proof of $\|R\mu\|^2 \gtrsim \sum_{j=0}^N\theta_j^2$} \label{secl22}

\subsection{The main lemma}

The main lemma to prove the estimate 
\begin{equation}\label{equp1}
\|R\mu\|^2 \gtrsim \sum_{j=0}^N\theta_j^2
\end{equation}
 is the following.

\begin{lemma}\label{mainlem}
We have
\begin{equation}\label{equp2}
\sum_{Q\in \Delta^0} \|D_Q (R\mu)\|^2\gtrsim \sum_{j=0}^{N-1}\theta_j^2.
\end{equation}
\end{lemma}

\medskip
Let us see how one deduces \rf{equp1} from the preceding inequality.

\bigskip
\begin{proof}[\bf Proof of \rf{equp1} using Lemma \ref{mainlem}]
From \rf{eqort1} and \rf{equp2} we infer that
\begin{equation}\label{eqpas1}
\|R\mu\|^2 \geq \|S_N(R\mu)\|^2\geq C_{3}^{-1} \sum_{j=0}^{N-1}\theta_j^2.
\end{equation}
So we only have to show that $\|R\mu\|^2 \gtrsim \theta_N^2.$

Consider $Q\in\Delta_N$ and $x\in Q$. We split $R\mu(x)$ as follows:
\begin{align*}
R\mu(x) & = R(\chi_Q\mu)(x) + R(\chi_{\R^d\setminus Q}\mu)(x)\\
& =
R(\chi_Q\mu)(x) + S_N(R\mu)(x) + \bigl(R(\chi_{\R^d\setminus Q}\mu)(x)
 - S_N(R\mu)(x)\bigr).
 \end{align*}
So we get
\begin{align}\label{eqd43}
\|R\mu\| &\geq \Bigl\|\sum_{Q\in\Delta_N}\chi_Q R(\chi_Q\mu)\Bigr\| -\|S_N(R\mu)\|
\nonumber\\
&\quad
-\Bigl\|\sum_{Q\in\Delta_N}\chi_Q R(\chi_{\R^d\setminus Q}\mu)
 - S_N(R\mu)\Bigr\|.
 \end{align}
It is easy to check that 
$$\Bigl\|\sum_{Q\in\Delta_N}\chi_Q R(\chi_Q\mu)\Bigr\|\geq C_4^{-1}\theta_N.$$ 
To deal with $S_N(R\mu)$ we simply use the fact that
$$\|S_N(R\mu)\|\leq \|R\mu\|.$$
On the other hand, by Lemma \ref{lemnab}, if $x\in Q\in\Delta_N$,
$$|R(\chi_{\R^d\setminus Q}\mu)(x)
 - S_N(R\mu)(x)| = |R(\chi_{\R^d\setminus Q}\mu)(x)
 - S_Q(R(\chi_{\R^d\setminus Q}\mu)(x)| \lesssim p_{N-1}.$$
By Cauchy-Schwartz, it follows easily that
$p_{N-1}\leq C\bigl(\sum_{j=0}^{N-1}\theta_j^2\bigr)^{1/2}$. Then we deduce
$$\Bigl\|\sum_{Q\in\Delta_N}\chi_Q R(\chi_{\R^d\setminus Q}\mu)
 - S_N(R\mu)\Bigr\|^2\leq 
C\sum_{j=0}^{N-1}\theta_j^2.$$
Then, by \rf{eqd43} and the estimates above, we get
$$\|R\mu\|\geq C_{4}^{-1}\theta_N - \|R\mu\|- C_{5}\bigl(\sum_{j=0}^{N-1}\theta_j^2\bigr)^{1/2}.$$
From \rf{eqpas1}, we infer that
$$C_{4}^{-1}\theta_N\leq  2\|R\mu\|+ C_{5}\Bigl(\sum_{j=0}^{N-1}\theta_j^2\Bigr)^{1/2}
\leq 2\|R\mu\|+ C_{3}^{1/2}\,C_{5}\|R\mu\|,$$
and thus the lemma follows.
\end{proof}


\subsection{The stopping scales and the intervals $I_k$}

To prove Lemma \ref{mainlem} we need to define some stopping scales on the squares from $\Delta$. Let $B$ be some big constant (say, $B>100$) to be fixed below.
We proceed by induction to define a subset $\sss:=\{s_0,\ldots,s_m\}\subset\{0,1,\ldots, N\}$. 
First we set  $s_0=0$.
If, for some $k\geq 0$, $s_{k}$ has already been defined and $s_{k}<N-1$, then $s_{k+1}$ is the least integer $i>s_k$ 
which verifies at least one of the following conditions:
\begin{itemize}
\item[(a)] $i=N$, or
\item[(b)] $\theta_i>B \,\theta_{s_{k}}$, or
\item[(c)] $\theta_i<B^{-1} \,\theta_{s_k}$.
\end{itemize}
We finish the construction of $\sss$ when we find some $s_{k+1}=N$. Notice that we have
$$[0,N-1]\cap \Z = \bigcup_{k=0}^{m-1} [s_k,s_{k+1})\cap \Z =: \bigcup_{k=0}^{m-1}
 I_k.$$
Moreover, the intervals $I_k$ are pairwise disjoint. 
 
If $s_k$ satisfies the condition (a) above, then we say that $I_k$ is terminal (in this case $k+1=m$). 
If $s_k$ satisfies (b) but not (a), then we say that
$I_k$ is an interval of increasing density, $I_k\in ID$. If (c) holds for $s_k$, but not (a) nor (b), then we say that
$I_k$ is  an interval of decreasing density, $I_k\in DD$. We denote its length by $|I_k|$. Notice that it coincides with $\# I_k$.

For $0\leq k\leq m$, we denote
$$T_k\mu = \sum_{j:s_k\leq j< s_{k+1}} D_j(R\mu).$$
In this way,
$$S_N(R\mu) = \sum_{k=0}^{m-1} T_k\mu,$$
and since the functions $D_j(R\mu)$ are pairwise orthogonal,
$$\|S_N(R\mu)\|^2 = \sum_{k=0}^{m-1}\|T_k\mu\|^2.$$

To simplify notation, given $A\subset \{0,\ldots,N\}$, we denote
$$\sigma(A) := \sum_{j\in A}\theta_j^2.$$
So $\sigma$ can be thought as a measure on $\{0,\ldots,N\}$.


\subsection{Good and bad scales}

We say that $j\in\{0,N-1\}$ is a good scale, and we write $j\in \GZ$, if
$$p_j\leq 40\theta_j.$$
Otherwise, we say that $j$ is a bad scale and we write $j\in \BZ$.

\begin{lemma}\label{lembons0}
We have
$$\sigma(\BZ)\leq \frac1{10}\,\sigma([0,N-1]).$$
\end{lemma}

\begin{proof}
As in \rf{eqpjtj} (replacing $M$ by $N-1$),
$$\sum_{j=0}^{N-1} p_j^2  \leq
 4\sum_{k=0}^{N-1} \theta_k^2 = 4\,\sigma([0,N-1]).$$
Thus,
$$\sigma(\BZ)\leq \frac1{40} 
\sum_{j=0}^{N-1} p_j^2 \leq
\frac1{10}\,\sigma([0,N-1]).$$
\end{proof}


\subsection{Good and bad intervals $I_k$}

We also say that an interval $I_k$ is good if
$$\sigma(I_k\cap \GZ) \geq \frac1{10}\,\sigma(I_k).$$ 
Otherwise we say that it is bad.

\begin{lemma}\label{lemgoodint}
$$\sigma([0,N-1])\leq \frac98\sum_{k:\,I_k {\rm \;good}}
\sigma(I_k).$$
\end{lemma}

\begin{proof}
If $I_k$ is bad, then
$$\sigma(I_k\cap \BZ) \geq \frac9{10}\,\sigma(I_k).$$
Thus,
$$
\sum_{k:\,I_k {\rm \;bad}} \sigma(I_k)
\leq \frac{10}9\sigma(\BZ)\leq \frac{10}9\,\,\frac1{10}\,\sigma([0,N-1]) =
\frac19\,\sigma([0,N-1]).
$$
Therefore,
\begin{align*}
\sigma([0,N-1]) & = 
\sum_{k:\,I_k {\rm \;good}} \sigma(I_k)
+\sum_{k:\,I_k {\rm \;bad}} \sigma(I_k)\\
&\leq \sum_{k:\,I_k {\rm \;good}} \sigma(I_k) + \frac19 \sigma([0,N-1]), 
\end{align*}
 and so
$$\sigma([0,N-1]) \leq \frac98 
\sum_{k:\,I_k {\rm \;good}} \sigma(I_k).$$
\end{proof}


\subsection{Long and and short intervals $I_k$}\label{subsechigh}

Let $N_L$ be some (big) integer to be fixed below. 
We say that an interval $I_k$ is long if 
$$|I_k|=s_{k+1} - s_k \geq N_L.$$
Otherwise we say that $I_k$ is short.


\subsection{Estimates for long good intervals $I_k$. The key lemma}

\begin{lemma}\label{lemj0}
Let $I_k$ be good, and set $j_0 = \min(I_k\cap \GZ)$. 
 Then,
$$j_0-s_k \leq \frac{10B^4}{1+10B^4}\,(s_{k+1} - s_k).$$
\end{lemma}

\begin{proof}
We denote $\ell=s_{k+1} - s_k$ and $\lambda=j_0-s_k$. Then we have
$$\sigma(I_k\cap \GZ) \leq B^2\theta_{s_k}^2(\ell-\lambda),$$
and also
$$\sigma(I_k\cap\BZ)\geq B^{-2}\theta_{s_k}^2\,\lambda.$$
Since $I_k$ is good, we have $\sigma(I_k\cap \BZ)\leq 10\sigma(I_k\cap \GZ)$, and so we infer that
$$\lambda\leq 10B^4(\ell-\lambda),$$
and the lemma follows.
\end{proof}

\begin{lemma}\label{lemaux11}
Let $0\leq k\leq N-1$. There exists some absolute constant $C_{6}$ such that if
\begin{equation}\label{eqassum3}
\frac{\ell_k}{\ell_{k-1}}\,p_{k-1}\leq C_{6}\bigl(\theta_k + \theta_{k+1}+\ldots +
\theta_{k+h}\bigr),
\end{equation}
then
$$\sum_{j=k}^{k+h}\|D_j(R\mu)\|^2\geq C_7^{-1}2^{-hd}
\bigl(\theta_k + \theta_{k+1}+\ldots +
\theta_{k+h}\bigr)^2.$$
\end{lemma}

\begin{proof}
Denote $f= \sum_{j=k}^{k+h} D_j(R\mu)$.
Take $P\in \Delta_{k+h+1}$ and $Q\in\Delta_k$ containing $P$. 
Then, for $x\in P$ we have
$$f(x) = S_P(R\mu)(x) - S_Q(R\mu)(x).$$
By antisymmetry, as in \rf{eqdec0}, we get
$$f(x) = S_P(R(\chi_{Q\setminus P}\mu))(x) + 
S_P(R(\chi_{\R^d\setminus Q}\mu))(x) - S_Q(R(\chi_{\R^d\setminus Q}\mu))(x).$$
From Lemma \ref{lemnab} it follows that
$$|S_P(R(\chi_{\R^d\setminus Q}\mu))(x) - S_Q(R(\chi_{\R^d\setminus Q}\mu))(x)|
\leq C_8\,\frac{\ell_k}{\ell_{k-1}}\,p_{k-1}.$$
On the other hand, if $P\in\Delta_{k+h+1}$ is a cube containing a corner of $Q$, then it is easy to
check that
$$\biggl|\frac1{\mu(P)} \int_P R(\chi_{Q\setminus P}\mu)\,d\mu\biggr| \geq C_9^{-1}
\bigl(\theta_k + \theta_{k+1}+\ldots
\theta_{k+h}\bigr).$$
Therefore,
$$|f(x)| \geq C_9^{-1}
\bigl(\theta_k + \theta_{k+1}+\ldots
\theta_{k+h}\bigr) - C_8\frac{\ell_k}{\ell_{k-1}}\,p_{k-1}.$$
As a consequence, if $C_{6}\leq C_9^{-1}C_8^{-1}/2$, then
\begin{align*}
\|\chi_Q \,f\|^2& \geq C^{-1}
\bigl(\theta_k + \theta_{k+1}+\ldots
\theta_{k+h}\bigr)^2\mu(P) \\
& = 2^{-(h+1)d}C^{-1}
\bigl(\theta_k + \theta_{k+1}+\ldots
\theta_{k+h}\bigr)^2\mu(Q).
\end{align*}
Summing over all the cubes $Q\in \Delta_k$, the lemma follows.
\end{proof}


\begin{lemma}\label{lemaux00}[Key lemma]
Let $A,c_0$ be positive constants, and $r,q\in[0,N-1]\cap\Z$ such that $q\leq r$,
$\,\frac{\ell_q}{\ell_{q-1}}\,p_{q-1}\leq c_0\theta_q$ and,
for all $j$ with $q\leq j\leq r$,
$$A^{-1}\theta_q\leq \theta_j\leq A\theta_q.$$
There exists $N_1=N_1(c_0,A)$ such that if $|q-r|>N_1$, then
$$\sum_{j=q}^r\|D_j(R\mu)\|^2\geq C|q-r|\,\theta_q^2,$$
where $C$ is some positive constant depending on $c_0$ and $A$.
\end{lemma}

\begin{proof}
Set $f = \sum_{j=q}^r D_j(R\mu)$. We have to show that $\|f\|^2
\geq C|q-r|\,\theta_q^2$. 

Let $M_0$ some positive integer depending on $c_0,A$ to be  fixed below.  We decompose $f$ as follows
\begin{equation}\label{eqsum83}
f =  \sum_{j=q}^{q+t\,M_0 -1}
D_j(R\mu) + \sum_{j=q+t\,M_0}^{r}D_j(R\mu),
\end{equation}
where $t$ is the biggest integer such that $q+t\,M_0 -1\leq r$.
Assuming $N_1$ big enough we will have $t\approx |q-r|$, with constants 
depending on $M_0$, and so on $c_0,A$.

We write the first sum on the right side of \rf{eqsum83} as follows:
$$ \sum_{j=q}^{q+t\,M_0 -1}
D_j(R\mu) = \sum_{h=0}^{t-1} \,\sum_{j=q+hM_0}^{q+(h+1)M_0 -1} D_j(R\mu)
=: \sum_{h=0}^{t-1} U_h(\mu).$$
By orthogonality, we have 
$$\|f\|^2 \geq \sum_{h=0}^{t-1} \|U_h(\mu)\|^2.$$
We will show below that if the parameter $M_0=M_0(c_0,A)$ is chosen big enough, then
\begin{equation}\label{eqcla94}
\|U_h(\mu)\|^2\geq C(c_0,A) \theta_{q}^2\quad \mbox{ for all $0\leq h\leq t-1$},
\end{equation}
and thus
$$\|f\|^2 \geq C(c_0,A) \,|q-r|\,\theta_q^2,$$
if $N_1\geq 2M_0$, say.

To prove \rf{eqcla94} we intend to apply Lemma \ref{lemaux11}. 
Recall that $\frac{\ell_q}{\ell_{q-1}}\,p_{q-1}\leq c_0 \theta_{q}$, and since
\begin{align*}
p_{q+hM_0-1} &= \sum_{i\leq q+hM_0-1} \frac{\ell_{q+hM_0-1}}{\ell_i}\,\theta_i\\
&=
\sum_{q-1<i\leq q+hM_0-1} \frac{\ell_{q+hM_0-1}}{\ell_i}\,\theta_i+ \frac{\ell_{q+hM_0-1}}{\ell_{q-1}}\,p_{q-1},
\end{align*}
we infer that
$$p_{q+hM_0-1}\leq 2A\theta_{q} + \frac{\ell_{q+hM_0-1}}{\ell_{q-1}}\,p_{q-1}.$$
Therefore, 
$$\frac{\ell_{q+hM_0}}{\ell_{q+hM_0-1}}\,p_{q+hM_0-1}
\leq 2A\theta_{q} + \frac{\ell_{q+hM_0}}{\ell_{q-1}}\,p_{q-1}
\leq 2A\theta_{q} + \frac{\ell_{q}}{\ell_{q-1}}\,p_{q-1}
\leq (2A+c_0)\theta_q.$$
On the other hand, 
$$\sum_{j=q+hM_0}^{q+(h+1)M_0 -1} \theta_j \geq M_0A^{-1}\theta_q.$$
If $M_0$ is big enough then $2A+c_0\leq C_6M_0A^{-1}$ and so the assumption 
\rf{eqassum3} in Lemma 
\ref{lemaux11} is satisfied. Thus
$$\|U_h\mu\|^2\geq C_7^{-1}2^{-M_0d} 
\biggl(\sum_{j=q+hM_0}^{q+(h+1)M_0 -1} \theta_j\biggr)^2\geq C_7^{-1}2^{-M_0d}M_0^2A^{-2}\theta_q^2,$$ 
and so our claim \rf{eqcla94} follows.
\end{proof}

\begin{lemma}\label{lemlongood}
Suppose that the constant $N_L$ is chosen big enough (depending on $B$).
If $I_k$ is long and good, then 
$$\sigma(I_k)\leq C(B)\|T_k\mu\|^2.$$
\end{lemma}

Recall that $T_k\mu = \sum_{j:s_k\leq j< s_{k+1}} D_j(R\mu)$. 

\begin{proof}
Set $\ell=s_{k+1} - s_k$. Notice that 
$$\sigma(I_k)\leq \ell\, B^2\theta_{s_k}^2.$$
Let $j_0 = \min(I_k\cap \GZ)$. 
We suppose that $N_L\gg B^4$, so that by Lemma \ref{lemj0},
$$s_{k+1}-j_0\geq \frac1{10B^4}\,\ell\gg 1.$$

  We split $T_k\mu$ as follows
\begin{equation*}
T_k\mu = \sum_{j=s_k}^{j_0-1} D_j(R\mu) + \sum_{j=j_0}^{s_{k+1}-1}
D_j(R\mu),
\end{equation*}
Now we apply Lemma \ref{lemaux00}, with $A=B$, $c_0=40$, and we we deduce that if 
$N_L$ is big enough, then
$$\sum_{j=j_0}^{s_{k+1}-1}
\|D_j(R\mu)\|^2\geq C(B)^{-1}|s_{k+1}-j_0|\,\theta_{s_k}^2$$
By orthogonality, 
$$\|T_k\mu\|^2\geq \sum_{j=j_0}^{s_{k+1}-1}
\|D_j(R\mu)\|^2,$$
and thus the lemma follows.
\end{proof}


\subsection{The intervals $J_h$}

By Lemmas \ref{lemgoodint} and \ref{lemlongood}, 
to finish our proof of $\sigma([0,N-1])\lesssim \sum_j\|D_j(R\mu)\|^2$, it is enough to show that
\begin{equation}\label{eqsho1}
\sum_{k:\,I_k {\rm \;short\; good}} \sigma(I_k) \lesssim \sum_j\|D_j(R\mu)\|^2.
\end{equation}
To this end, we have to define some auxiliary intervals $J_h$.

We consider the following partial ordering in the family of intervals contained in $\R$:
if $I,J$ are disjoint intervals such that all $x\in I$, $y\in J$ satisfy $x<y$, then we write $I\prec J$.

An interval $J_h$, $h\geq 1$, is the union of two intervals $I_k,I_{k+1}$, so that $I_k$ is of type
$ID$ and $I_{k+1}$ is either of type $DD$ or it is the terminal interval $I_m$. 
Then $\{J_h\}_{1\leq h\leq m_J}$ is the collection
of all these intervals. We assume that $J_h\prec J_{h+1}$ for all $h$.
Moreover, for convenience, if $I_0$ is of type $DD$, we set $J_0=I_0$.

\begin{rem}\label{rem58}
Of course, there may be intervals $I_k$ which are not contained in any interval $J_h$. 
Suppose that, for some $0\leq h\leq m_J$, there are intervals $I_k$ such that
$$J_h\prec I_k\prec I_{k+1}\prec\ldots \prec I_{k+r}\prec J_{h+1}.$$
Then, from the definition of the intervals $J_h$, it turns out that either all the intervals $I_{k},\ldots,I_{k+r}$ are of type $ID$, or all are of type $DD$, or 
there exists 
$1\leq s\leq r$ such that $I_{k},\ldots,I_{k+s-1}$ are of type $DD$, and 
$I_{k+s},\ldots,I_{k+r}$ are of type $ID$.
\end{rem}

Given an interval $I\subset[0,N]$, we denote
$$\theta^{\max}(I)=\max_{j\in I}\theta_j.$$

\begin{lemma}\label{lemamax11}
Let $J_h$, $0\leq h\leq m_J-1$, be such that
$$J_h\prec I_k\prec I_{k+1}\prec\ldots \prec I_{k+r}\prec J_{h+1},$$
or in the case $h=m_J$,
$$J_h\prec I_k\prec I_{k+1}\prec\ldots \prec I_{k+r}.$$
Then,
\begin{equation}\label{eqkid}
\sum_{\substack{k\leq i\leq k+r \\ I_i{\; \rm short}}} \sigma(I_i)
 \leq C(B,N_L) \bigl[\theta^{\max}(J_h)^2 + \theta^{\max}(J_{h+1})^2\bigr],
\end{equation}
where, for convenience, we set $\theta_{m_J+1}=0$.
\end{lemma}

\begin{proof}
 Notice that any short interval $I_k$
 satisfies
\begin{equation}\label{eqddid}
\sigma(I_k) \leq B^2\,N_L\theta_{s_k}^2.
\end{equation}
If there is some $q\geq1$ such that the intervals $I_k,\ldots,I_{k+q-1}$ are of type
$DD$, then
$$\theta_{s_{k+q-1}}\leq B^{-1}
\theta_{s_{k +q-2}}\leq\ldots \leq B^{1-q}\theta_{s_{k}}\leq B^{-q}\theta^{\max}(J_h).$$
Thus,
$$\sum_{\substack{k\leq i\leq k+q-1 \\ I_i{\; \rm short}}} \sigma(I_i)
 \leq C(B,N_L)\theta^{\max}(J_h)^2.$$
 Analogously, one deduces that
$$\sum_{\substack{k+q\leq i\leq k+r \\ I_i{\; \rm short}}} \sigma(I_i)
 \leq C(B,N_L)\theta^{\max}(J_{h+1})^2,$$
and the lemma follows.
\end{proof}

\begin{lemma}\label{lemjh}
We have
\begin{equation}\label{eqkid2}
\sum_{\substack{k:I_k{\; \rm short}}} \sigma(I_k)
 \leq C(B,N_L)\sum_{h=0}^{m_J} \theta^{\max}(J_h)^2.
\end{equation}
\end{lemma}

\begin{proof}This is a direct consequence of Lemma \ref{lemamax11}.
\end{proof}


\subsection{The standard intervals $J_h$}

By Lemma \ref{lemjh}, in order to prove \rf{eqsho1}, it is enough to show that
\begin{equation*}
\sum_{h} \theta^{\max}(J_h)^2 \lesssim \sum_j\|D_j(R\mu)\|^2.
\end{equation*} 
To this end, we need to 
distinguish different types of intervals $J_h$. 
For $h\geq1$, let $t_h\in J_h$ be the least integer such that
$$\theta_{t_h}>B^{-1/2}\,\theta^{\max}(J_h).$$
Notice that, if $J_h= I_k\cup I_{k+1}$, then
$\theta^{\max}(J_h)\leq B\theta_{s_{k+1}}$. However we cannot
ensure that $\theta^{\max}(J_h)\leq B^2\theta_{s_{k}}$ because
it may happen that $\theta_{s_{k+1}}\gg B\theta_{s_k}$.

We say that $J_h$ is {\bf standard} if 
 \begin{equation}\label{eqpekt}
\frac{\ell_{t_h}}{\ell_{t_h-1}}\,p_{t_h-1} \leq
C_{10}\,\theta^{\max}(J_h),
\end{equation}
 where $C_{10}=C_6/2$ (with $C_6$ from \rf{eqassum3}.
For convenience, if $J_0$ exists (and thus $J_0=I_0\in DD$) we also say that $J_0$ is standard.

\begin{lemma}\label{lemstan}
If $J_h$ is standard, then
$$\theta^{\max}(J_h)^2\leq  C(B) \sum_{j\in J_h} \|D_j(R\mu)\|^2.$$
\end{lemma}

\begin{proof}
In the special case $h=0$ (with $J_0=I_0$), it is immediate to check that $\|D_0(R\mu)\|^2\geq C^{-1}\theta_0^2 \geq C^{-1}B^{-2}
 \theta^{\max}(J_0)^2$ (for instance, one can apply Lemma \ref{lemaux11} with $p_{-1}=0$), 
 and thus the lemma holds.

For $h\geq 1$, we set
$$J_h = [s_k,t_h-1) \cup[t_h,s_{k+2}) =: J_h^a \cup J_h^b.$$
Observe that $\theta_{max}(J_h)$ is attained at some scale from $J_h^b$, and $\theta_j\leq B^{-1/2}
\theta_{max}(J_h)$ for $j\in J_h^a$.


We distinguish two cases: 

\vv
\noi {\bf Case 1.} Suppose first that the length $|J_h^b|$ is big. That is,  $|J_h^b|=s_{k+2} - t_h>N_2$, where $N_2=N_2(C_{10},B)$ is some big integer.
By \rf{eqpekt}, we have
$$\frac{\ell_{t_h}}{\ell_{t_h-1}}\,p_{t_h-1} \leq
C_{10}\,\theta^{\max}(J_h)\leq C(B) \theta_{t_h},$$
 and thus from Lemma 
\ref{lemaux00} we infer that if $N_2$ is chosen big enough, then
$$\theta^{\max}(J_h)^2\leq  C(B)\sum_{j\in J_h^b} \|D_j(R\mu)\|^2,$$
and so the lemma holds in this case.

\vv
\noi {\bf Case 2.} 
Assume that $|J_h^b|\leq N_2$. 
From \rf{eqpekt}, recalling that $C_{10}=C_6/2$, we infer that
$$\frac{\ell_{t_h}}{\ell_{t_h-1}}\,p_{t_h-1} \leq C_{10}\,\theta^{\max}(J_h)\leq
C_6\sum_{j\in J_h^b} \theta_j,$$
and then, by Lemma \ref{lemaux11},
$$\sum_{j\in J_j^b}\|D_j(R\mu)\|^2\geq C_7^{-1}2^{-N_2d}
\Bigl(\sum_{j\in J_h^b}\theta_j\Bigr)^2
\geq C_7^{-1}2^{-N_2d}\,\theta^{\max}(J_h)^2,$$
and so the lemma also holds in this situation.
\end{proof}


\subsection{The non standard intervals $J_h$}

\begin{lemma}\label{lemnonstan}
Suppose that $B$ has been chosen big enough. We have
$$\sum_{h:J_h{\; \rm non\; standard}}\theta^{\max}(J_h)^2\leq  C(B) 
\sum_{h:J_h{\; \rm standard}}\theta^{\max}(J_h)^2.$$
\end{lemma}

\begin{proof}
Denote by $\{J^{st}_n\}_n$ the subfamily of the standard intervals from $\{J_h\}_h$, ordered so that
$J^{st}_n\prec J^{st}_{n+1}$ for all $n$. For a fixed $n$, denote by $\Lambda_1$, \ldots, $\Lambda_m$
the collection of all non standard intervals from the family $\{J_h\}$ 
such that either
$$J^{st}_n\prec\Lambda_1\prec\Lambda_2\prec\ldots\prec \Lambda_m\prec J^{st}_{n+1} \qquad\mbox{if $J^{st}_{n+1}$ exists,}$$
or
$$J^{st}_n\prec\Lambda_1\prec\Lambda_2\prec\ldots\prec \Lambda_m\qquad\mbox{if $J^{st}_{n+1}$ does not exist.}$$
We will prove that
\begin{equation}\label{eqclau23}
\theta^{\max}(\Lambda_i)\leq B^{-i/8s} \,\theta^{\max}(J^{st}_n)
\quad\mbox{ for $i\geq1$},
\end{equation}
by induction on $i$. The lemma follows easily from this estimate.

To simplify notation, we set $\Lambda_0=J^{st}_n$ and $\theta^{\max}_i=\theta^{\max}(\Lambda_i)$.
Also, if $\Lambda_i=I_k\cup I_{k+1}$, we denote by $Q_i$ a cube from $\Delta_{s_k}$, by $\wt Q_i$ 
a cube from $\Delta_{t_h-1}$ (see \rf{eqpekt}), 
and by $Q_i^{\max}$ a cube from $\bigcup_{j\in I_{k+1}}\Delta_j$ such that $\theta^{\max}_i=\theta_j$. Moreover, we assume that $$Q_i\supset \wt Q_i\supset Q_i^{\max}\supset Q_{i+1}
\supset
 \wt Q_{i+1}\supset Q_{i+1}^{\max}\supset\ldots$$

First
we prove \rf{eqclau23} for {\boldmath $i=1$}. Since $\Lambda_1$ is
not standard,
\begin{equation}\label{eqnos3}
\theta_1^{\max}  \leq C_{10}^{-1}\,\frac{\ell(s(\wt Q_1))}{\ell(\wt
Q_1)}\,p(\wt Q_1),
\end{equation}
where $s(\wt Q_1)$ stands for a son of $\wt Q_1$.
To estimate $p(\wt Q_1)$ (recall the notation in \rf{defpqr}), we decompose it as follows:
$$p(\wt Q_1) \leq p(\wt Q_1,Q_1) + \frac{\ell(\wt Q_1)}{\ell(Q_1)}\,p(Q_1,Q_0^{\max}) + \frac{\ell(\wt Q_1)}{\ell(Q_0^{\max})}\,p(Q_0^{\max}).$$
Now observe that
\begin{equation}\label{eqpq1q1}
p(\wt Q_1,Q_1) \leq 2B^{-1/2}\theta_1^{\max},
\end{equation}
since $\theta(P)\leq B^{-1/2}\theta_1^{\max}$ for $\wt Q_1\subset
P\subset Q_1$. Also,
\begin{equation}\label{eqpq1q0}
p(Q_1,Q_0^{\max})\leq 2\theta(Q_1) + 2\theta_0^{\max},
\end{equation}
 because
$\theta(P)\leq \theta(Q_1) + \theta_0^{\max}$ for $Q_1\subset
P\subset Q_0^{\max}$, taking into account Remark \ref{rem58}.
 And finally,
\begin{align}\label{eqpk00}
p(Q_0^{\max})& \leq p(Q_0^{\max},\wt Q_0) +\frac{\ell(Q_0^{\max})}{\ell(\wt Q_0)}\,p(\wt Q_0) \\
& \leq p(Q_0^{\max},\wt Q_0) +\frac{\ell(s(\wt Q_0))}{\ell(\wt
Q_0)}\,p(\wt Q_0) \leq 4\theta_0^{\max},\nonumber
\end{align}
because $\theta(P)\leq \theta_0^{\max}$ for $Q_0^{\max}\subset P\subset \wt Q_0$ and moreover $\Lambda_0$
is standard (we assume $C_6\leq 1$, say). Thus we infer that
\begin{align*}
p(\wt Q_1) & \leq 2B^{-1/2}\theta_1^{\max} + \frac{2\ell(\wt Q_1)}{\ell(Q_1)}\,\bigl(\theta(Q_1) +
\theta_0^{\max}\bigr) + \frac{4\ell(\wt Q_1)}{\ell(Q_0^{\max})}\,\theta_0^{\max}\\
& \leq 4B^{-1/2}\theta_1^{\max} +  \frac{6\ell(\wt
Q_1)}{\ell(Q_1)}\,\theta_0^{\max},
\end{align*}
using that $\theta(Q_1)\leq B^{-1}\theta_1^{\max}\leq B^{-1/2}\theta_1^{\max}$ in the second inequality.
If we plug this estimate into \rf{eqnos3} we deduce
$$
\theta_1^{\max}  \leq 4C_{10}^{-1}B^{-1/2}\theta_1^{\max} +
6C_{10}^{-1} \frac{\ell(s(\wt Q_1))}{\ell(Q_1)}\,\theta_0^{\max}.
$$
If we assume $B$ big enough, so that $4C_{10}^{-1}B^{-1/2}\leq
1/2$ (recall that $C_{10}=C_6/2$ does not depend on $B$), we obtain
$$
\theta_1^{\max}  \leq 12C_{10}^{-1} \frac{\ell(s(\wt
Q_1))}{\ell(Q_1)}\,\theta_0^{\max}.
$$
On the other hand, since $\theta(s(\wt Q_1))>B^{1/2}\theta(Q_1)$ (by the definition of $\wt Q_1$), we
infer that
\begin{equation}\label{eqll2}
\ell(s(\wt Q_1))^s\leq B^{-1/2}\ell(Q_1)^s,
\end{equation}
 and so
$$
\theta_1^{\max}  \leq 12C_{10}^{-1}B^{-1/2s}\,\theta_0^{\max}.
$$
If we suppose $B$ big enough again, \rf{eqclau23} follows in the particular case $i=1$.

The proof of \rf{eqclau23} for an arbitrary integer {\boldmath
$i\geq2$ }when we assume that it holds for $1,\ldots,i-1$ is
analogous to the one for the case $i=1$. For the sake of completeness
we will show the detailed arguments. As in \rf{eqnos3}, we
have
\begin{equation}\label{eqnosi}
\theta_i^{\max}  \leq C_{10}^{-1}\frac{\ell(s(\wt Q_i))}{\ell(\wt
Q_i)}\,p(\wt Q_i),
\end{equation}
because $\Lambda_i$ is not standard. Now we split $p(\wt Q_i)$ as follows:
\begin{align*}
p(\wt Q_i) & \leq p(\wt Q_i,Q_i) + \frac{\ell(\wt
Q_i)}{\ell(Q_i)}\,p(Q_i,Q_{i-1}^{\max}) \\
& \quad + \sum_{j=1}^{i-1}\frac{\ell(\wt
Q_i)}{\ell(Q_j^{\max})}\,p(Q_j^{\max},Q_{j-1}^{\max}) +
\frac{\ell(\wt Q_i)}{\ell(Q_0^{\max})}\,p(Q_0^{\max}).
\end{align*}
 We will
estimate each of the terms in the preceding inequality separately.
As in \rf{eqpq1q1}, we have
 $$p(\wt Q_i,Q_i) \leq 2B^{-1/2}\theta_i^{\max},$$
and as in \rf{eqpq1q0},
$$p(Q_i,Q_{i-1}^{\max})\leq 2\theta(Q_i) +
2\theta_{i-1}^{\max} \leq 2B^{-1/2}\theta_i^{\max} +
2\theta_{i-1}^{\max}.$$ By analogous arguments,
$$p(Q_j^{\max},Q_{j-1}^{\max}) \leq 2\theta_j^{\max} +
2\theta_{j-1}^{\max}.$$ On the other hand, the term
$p(Q_0^{\max})$ has been estimated in \rf{eqpk00}. By the
preceding inequalities and the induction hypothesis, we obtain
\begin{align*}
p(\wt Q_i) & \leq 2B^{-1/2}\theta_i^{\max} + \frac{\ell(\wt
Q_i)}{\ell(Q_i)}\,\bigl( 2B^{-1/2}\theta_i^{\max} +
2\theta_{i-1}^{\max}\bigr) \\
&\quad + \sum_{j=1}^{i-1}\frac{\ell(\wt
Q_i)}{\ell(Q_j^{\max})}\,\bigl(2\theta_j^{\max} +
2\theta_{j-1}^{\max}\bigr) +\frac{4\ell(\wt
Q_i)}{\ell(Q_0^{\max})}\theta_0^{\max} \\
& \leq 4B^{-1/2}\theta_i^{\max} + 2\,\frac{\ell(\wt
Q_i)}{\ell(Q_i)}\,B^{-(i-1)/8s}\,\theta_0^{\max}\\
&\quad + 4\sum_{j=1}^{i-1}\frac{\ell(\wt
Q_i)}{\ell(Q_j^{\max})}\,B^{-(j-1)/8s}\,\theta_0^{\max} +
\frac{4\ell(\wt Q_i)}{\ell(Q_0^{\max})}\theta_0^{\max}.
\end{align*}
If we plug this inequality into \rf{eqnosi} and we assume $B$ big
enough, we deduce that
\begin{align} \label{eqll3}
\theta_i^{\max}  & \leq C\, \biggl[ \frac{\ell(s(\wt
Q_i))}{\ell(Q_i)}\,B^{-(i-1)/8s}\,\theta_0^{\max}\\
&\quad + \sum_{j=1}^{i-1}\frac{\ell(s(\wt
Q_i))}{\ell(Q_j^{\max})}\,B^{-(j-1)/8s}\,\theta_0^{\max} +
\frac{\ell(s(\wt Q_i))}{\ell(Q_0^{\max})}\,\theta_0^{\max}\biggr],
\nonumber
\end{align}
with $C$ independent of $B$.
As in \rf{eqll2}, we have
$$\frac{\ell(s(\wt Q_i))}{\ell(Q_i)} \leq B^{-1/2s},$$
and for $0\leq j\leq i-1$,
$$\frac{\ell(s(\wt Q_i))}{\ell(Q_j^{\max})} \leq
\frac{\ell(s(\wt Q_i))}{\ell(Q_i)}\cdots\frac{\ell(s(\wt
Q_{j+1}))}{\ell(Q_{j+1})} \leq B^{(j-i)/2s}.
$$
From the latter estimates and \rf{eqll3} we obtain
\begin{align*}
\theta_i^{\max}  & \leq C\, \biggl[ B^{-1/2s}\,B^{-(i-1)/8s}\,\theta_0^{\max} \\
&\quad + \sum_{j=1}^{i-1}
B^{(j-i)/2s}\,B^{-(j-1)/8s}\,\theta_0^{\max}
+  B^{-i/2s}\,\theta_0^{\max}\biggr] \\
& \leq C\,B^{-1/4s}\,B^{-i/8s}\,\theta_0^{\max},
\end{align*}
and so \rf{eqclau23} holds if we assume $B$ big enough. 
\end{proof}


\subsection{Proof of Lemma \ref{mainlem}}

From Lemmas \ref{lemgoodint}, \ref{lemlongood}, and \ref{lemjh}, we get
\begin{align*}
\sum_{j=0}^{N-1}\theta_j^2 & = \sum_k \sigma(I_k) \leq
C\sum_{k:\,I_k{\rm \; good}} \sigma(I_k)\\
& =
C\sum_{k:\,I_k{\rm \;long\; good}} \sigma(I_k) + 
C\sum_{k:\,I_k{\rm \;short\; good}} \sigma(I_k)\\
&\leq C\sum_{j=0}^{N-1}\|D_j(R\mu)\|^2 + C\sum_h \theta^{\max}(J_h)^2.
\end{align*}

By Lemmas \ref{lemnonstan} and \ref{lemstan},
$$\sum_{h:\,J_h}
\theta^{\max}(J_h)^2\lesssim\sum_{h:\,J_h{\rm \;standard}}
\theta^{\max}(J_h)^2\lesssim\sum_{j=0}^{N-1} \|D_j(R\mu)\|^2.$$
We are done.
\fiproof


\section{Open problems}\label{secopen}

In this section we discuss some open problems in connection with 
Riesz transforms and Wolff potentials.


\vv
\vv 
\noi {\bf 1) Riesz transforms and rectifiability.}

Let $E\subset\R^d$ be a compact set with $0<\HH^n(E)<\infty$, for some integer $0<n<d$, and set $\mu=\HH^n_{|E}$. If $R_\mu^n$ is bounded in $L^2(\mu)$, is then $E$ $n$-rectifiable?
Recall that $E$ is called $n$-rectifiable if there exist Lipschitz mappings $g_i:\R^n\to\R^d$ such that 
$$\mu\Bigl(\R^d\setminus \bigcup_{i=1}^\infty g_i(\R^n)\Bigr) = 0.$$
When $n=1$, David and L\'eger \cite{Leger} answered the question in the affirmative, using the relationship between curvature and the Cauchy kernel.
By \cite{Volberg}, when $n=d-1$ this question is equivalent
to the following: is it true that $\kappa(E)=0$ if and only if $E$ is purely $(d-1)$-unrectifiable? ($E$ is called
purely $(d-1)$-unrectifiable if it does not contain any $n$-rectifiable subset $F$ with $\HH^{d-1}(F)>0$).

A partial result was obtained in \cite{Tolsa-vpriesz}, where it was shown that the existence of the principal values
$\lim_{\ve\to0} R_\ve^n\mu(x)$ for $\mu$-a.e. $x\in \R^d$ implies $E$ to be $n$-rectifiable. Under the additional
assumption
\begin{equation}\label{eqdens4}
\theta_{\mu,*}^n(x):= 
\liminf_{r\to0}\frac{\mu(B(x,r))}{r^n}>0\qquad \mbox{$\mu$-a.e. on $\R^d$},
\end{equation}
this had been proved previously by Mattila and Preiss in \cite{Mattila-Preiss}. Unfortunately, it is not known if
 the $L^2(\mu)$ boundedness of the Riesz transform $R^n_\mu$ implies the existence
of principal values, and so the results in \cite{Tolsa-vpriesz} and \cite{Mattila-Preiss} do not help to
solve the problem above.

Another related result is given in
\cite[Theorem 5.5]{Mattila-Preiss}, where it is proved  
that if \rf{eqdens4} holds and all the operators 
$$Tf(x) = \int K(x-y)f(y)\,d\mu(x),$$
with kernel of the form $K(x)=\vphi(|x|) x/|x|^{n+1}$ satisfying 
$|\nabla^j K(x)|\leq \frac{C(j)}{|x|^{n+j}}$ for $j\geq0$ 
 are bounded in 
$L^2(\mu)$, then $E$ is $n$-rectifiable.

A variant of this problem, posed by David and Semmes, consists in taking $E$ Ahlfors-David regular and $n$-dimensional. 
That is,
$$\HH^n(E\cap B(x,r))\approx r^n\quad\mbox{ for all $x\in E$, $0<r\leq\diam(E)$.}$$
Again, set $\mu=\HH^n_{|E}$.
If $R_\mu^n$ is bounded in $L^2(\mu)$, is then $E$ 
uniformly $n$-rectifiable? For the definition of uniform rectifiability, see \cite{David-Semmes-llibre1} and \cite{David-Semmes-llibre2} (for the reader's convenience let us say
that, roughly speaking, uniform rectifiability is the same as rectifiability plus some quantitative estimates). 
For $n=1$ the answer is true again, because of curvature. 
The result is from Mattila, Melnikov and Verdera \cite{Mattila-Melnikov-Verdera}. For $n> 1$, in \cite{David-Semmes-llibre1} and \cite{David-Semmes-llibre2} some partial answers 
are given. In particular, it is shown that if all the operators $T$ with kernel $K$ as above are bounded in $L^2(\mu)$,
then $E$ is uniformly rectifiable.


\vv
\vv 
\noi {\bf 2) 
Calder\'on-Zygmund capacities and Wolff potentials of non integer dimension.} 

This problem was already mentioned in the Introduction: is it true that for $0<s<d$ non integer we have
\begin{equation}\label{eqcomp3}
\gamma_s(E) \approx \dot C_{\frac23(d-s),\frac32}(E)
\end{equation}
with constants independent of $E$?
Recall that this was shown to be true when $0<s<1$ by Mateu, Prat and Verdera \cite{Mateu-Prat-Verdera}.  
In Theorem \ref{teokappa} we
have proved that \rf{eqcomp3}
holds for the particular case of the Cantor sets $E_N=E_N(\lambda)$ associated 
to a sequence $\lambda=(\lambda_n)^\infty_{n=1}$, with
$\lambda_n\le\tau_0<\frac{1}{2}.$ It might be also interesting to consider 
more general
Cantor type sets. For instance, for $n\geq0$ and $1\leq j\leq 2^{nd}$ take dilation factors $\lambda_{n,j}$, with $0<\lambda_{n,j}\leq \tau<2$,
and construct a Cantor type set analogous to the one in Theorem \ref{teokappa}, but
allowing different values $\lambda_{n,j}$, $1\leq j \leq2^{nd}$, for the different squares from a fixed scale $n$. Does \rf{eqcomp3} hold for these new Cantor sets?
To solve this more general case, one might try to implement the techniques 
used in Theorem \ref{teokappa}. However, several non trivial  difficulties arise.
First, one has to take into account that in this case, if one 
considers a probability measure $\mu$ such that $\mu(Q_{n,j})=2^{-nd}$ for each 
square $Q_{n,j}$ from the corresponding set $E_N$, then both estimates
$$\gamma_s(E_N) \approx \Bigl(\sum_{Q\in\Delta} \theta(Q_{n,j})^2\Bigr)^{-1/2},\qquad
\dot C_{\frac23(d-s),\frac32}(E_N) \approx \Bigl(\sum_{Q\in\Delta} \theta(Q_{n,j})^2\Bigr)^{-1/2}$$
{\bf are false} in general (as in Theorem \ref{teokappa}, $\Delta$ stands for the collection of all squares $Q_{n,j}$ in the construction of $E_N$). Thus $\mu$ should
be replaced by another measure. On the other hand, for some of the arguments involved in the estimates of the $L^2(\mu)$ norm of $R(\mu)$ in Theorem \ref{difi0},
the homogeneity of the set $E_N$ (i.e. the fact that $\ell(Q_{n,j})$ only depends on $n$) is essential.

As mentioned in the Introduction, it is proved in \cite{Ei-Na-Vo} that the estimate $\gamma_s(E) \gtrsim \dot C_{\frac23(d-s),\frac32}(E)$ holds for $0<s<d$.
The main obstacle to prove the opposite inequality is the following.
It is not known if, for $s\not\in\Z$, there are sets
$E$ with $0<\HH^s(E)<\infty$ such that the Riesz transform $R^s_\mu$, with
$\mu=\HH^s_{|E}$, is bounded in $L^2(\mu)$. If \rf{eqcomp3} holds, then such
sets do not exist. 
This is the case for $0<s<1$, as shown by Prat \cite{Prat-IMRN} 
using the curvature
method, and for other $s\not\in\Z$ by Vihtila \cite{Vihtila} 
under the additional assumption that
$\theta_{\mu,*}^s(x)>0$ for $\mu$-a.e. $x\in\R^d$, where $\theta_{\mu,*}^s(x)$
is defined in \rf{eqdens4}.

On the other hand, in \cite{Ruiz-Tolsa} it has been proved that, for $0<s<d$
and $\mu=\HH^s_{|E}$, with $0<\HH^s(E)<\infty$, the 
existence of the principal values
$\lim_{\ve\to0} R_\ve^s\mu(x)$ for $\mu$-a.e. $x\in \R^d$
 forces $s$ to be integer.
Notice that if one combines the results on principal values from 
\cite{Tolsa-vpriesz} mentioned above with the ones from \cite{Ruiz-Tolsa}, 
then one gets:

\begin{theorem*}\label{junt}
For $0< s \leq d$, let $E \subset \R^d$ be a set satisfying $0<\HH^s(E)<\infty$.
The principal value
$$\lim_{\ve \to 0} \int_{|x-y|>\ve}\frac{x-y}{|x-y|^{s+1}}d\HH^s_{|E}(y)$$ exists
for $\HH^s$-almost every $x \in E$ if and only if $s$ is integer and $E$ is $s$-rectifiable.
\end{theorem*}

It is interesting to compare the last theorem with well known results in geometric measure theory due essentially to Marstrand \cite{Marstrand} and Preiss \cite{Preiss}:

\medskip
\noindent{\em For $0< s \leq m$, let $E \subset \R^m$ be a set satisfying 
$0<\HH^s(E)<\infty$.
The density $\theta^s_{\HH^s|E}(x)$ exists 
for $\HH^s$-almost every $x \in E$ if and only if $s$ is integer and $E$ is $s$-rectifiable.}


\vv
\vv

\noi{\bf 3) $L^2$ boundedness of Riesz transforms and square functions.}

Given a non-increasing radial $\CC^\infty$ function $\psi$ such
that $\chi_{B(0,1/2)}\leq \psi\leq \chi_{B(0,2)}$, for each
$j\in\Z$, we set $\psi_j(z) := \psi(2^jz)$ and
$\vphi_j:=\psi_j-\psi_{j+1}$,  so that each function $\vphi_j$ is
non-negative and supported in the annulus $A(0,2^{-j-2},2^{-j+1})$, 
and moreover we have $\sum_{j\in\Z}\vphi_j(x) = 1$
for all $x\in\R^d\setminus \{0\}$. For
each $j\in Z$ we denote $K_j^s(x) = \vphi_j(x)\,x/|x|^{s+1}$ and
\begin{equation}\label{eqaux10}
R_j^s \mu(x) = \int K_j^s(x-y)\,d\mu(y).
\end{equation}
Notice that, at a formal level, we have $ R\mu = \sum_{j\in\Z} R_j\mu,$ and so
$$\|R^s\mu\|_{L^2(\mu)}^2 = \sum_{j\in \Z}\|R_j^s\mu\|_{L^2(\mu)}^2 + \sum_{j\neq k}\langle R_j^s\mu,\,R_k^s\mu\rangle.$$
Consider the square function
$$Q^s\mu(x) = \Bigl(\sum_{j\in\Z} |R_j^s\mu(x)|^2\Bigr)^{1/2},$$
and set $Q_\mu^s(f) = Q^s(f\,d\mu)$. Notice that 
$$\|Q^s_\mu(f)\|_{L^2(\mu)}^2= \sum_{j\in \Z}\|R_j^s(f\,d\mu)\|_{L^2(\mu)}^2.$$
One should view $Q^s_\mu(f)$ as a square function associated to the Riesz transform 
$R^s_\mu(f)$.

When $s$ is integer and $E\subset\R^d$ uniformly rectifiable, with $\mu= \HH^s_{|E}$, then
$Q^s_\mu$ is bounded in $L^2(\mu)$. Moreover, the converse is also true: 
if $E$ is Alhfors-David regular, the $L^2(\mu)$
boundedness of $Q_\mu$ implies that $E$ is uniformly rectifiable (at least for an appropriate
choice of the function $\psi$ above), as shown in \cite{Tolsa-alfas}. 
In the non Ahlfors-David regular case it is also true that the boundedness of $Q_\mu$ implies the rectifiability
of $E$ \cite{May-Volberg1}.

On the other hand, given
$E\subset\R^d$ such that $0<\HH^s(E)<\infty$, $0<s<d$, and $\mu= \HH^s_{|E}$, if $Q_\mu$
is bounded in $L^2(\mu)$, then $s\in\Z$. This follows easily from the results of
\cite{Ruiz-Tolsa}, as shown in \cite{May-Volberg2}. 
Thus the following question arises naturally:

\medskip
\noi {\em Let $0<s<d$ and let $\mu$ be a Radon measure on $\R^d$ with no atoms. 
Is it true that $R^s_\mu$ is bounded in 
$L^2(\mu)$ if and only if $Q^s_\mu$ is bounded in $L^2(\mu)$?}
\medskip

As remarked above, solving this question would be a fundamental contribution for the
solution of the problems explained above in 1) and 2).

\vv\vv
\noi{\bf 4) Bilipschitz and affine invariance, and other problems.}

Let $\mu$ be a Radon measure on $\C$ such that the Cauchy transform $\CC_\mu$
is bounded in $L^2(\mu)$. Recall that
$$\CC_\mu f(z)=\int\frac1{z-\xi}\,f(\xi)\,d\mu(\xi).$$
 In \cite{Tolsa-bilip} it has been shown that if
$\vphi:\C\to\C$ is a bilipschitz map and 
$\sigma=\vphi\#\mu$ is the image measure of $\mu$, then $\CC_\sigma$
is bounded in $L^2(\sigma)$. The analogous problem for the $(d-1)$-dimensional Riesz transform $R^{d-1}_\mu$ in $\R^{d}$ is open, and it seems that before trying to
solve it, one should understand better the relationship between the $L^2$ boundedness of the Riesz transforms and rectifiability [i.e.\ one should first solve
the questions in 1)], since this is a basic ingredient in the proof of 
the analogous result for the Cauchy transform in 
\cite{Tolsa-bilip}. However, in the case $d>2$, the problem is open even when 
$\vphi$ is an affine map. For instance, let 
$$\vphi(x_1,x_2,x_3,\ldots,x_d)= (2x_1,x_2,x_3,\ldots,x_d).$$
If $R^{d-1}_\mu$ is bounded in $L^2(\mu)$ and we set $\sigma=\vphi\#\sigma$, is then $R^{d-1}_\sigma$  bounded
in $L^2(\sigma)$?

A similar question in terms of the capacity $\kappa$ is the following. Is it
true that for any compact set $E\subset\R^d$, $\kappa(E)\approx\kappa(\vphi(E))$?
Analogous questions can be posed for the other capacities $\gamma_s$ and the Riesz transforms of codimension different from $1$.

\bigskip
Let us discuss another problem whose solution may help to understand
the relationship between the $L^2$ boundedness of Riesz transforms and geometry.
Let $R^s_{(j)}$, $0\leq j\leq d$, denote the scalar components of the (vectorial)
Riesz transform $R^s$. Let $\mu$ be a Radon measure on $\R^d$ such that
$\mu(B(x,r))\leq C\,r^s$ for all $x\in\R^d$, $r>0$. Suppose that $d-1$ components
of $R_\mu^s$, say $R_{(1),\mu}^s,\ldots,R_{(d-1),\mu}^s$ are bounded in $L^2(\mu)$.
Is then $R^s_\mu$ bounded in $L^2(\mu)$?
When $s=1$ the answer is yes, because of the curvature method. However, for other
values of $s$, the problem is open again. 

An analogous question can be posed in terms of
the capacities associated to these kernels. That is, for a compact set 
$E\subset\R^d$, let
$\wt \gamma_s(E) =  \sup|\nu(E)|$, where the supremum is taken over signed
measures (or distributions) supported on $E$ such that 
$\|R^s_{(j)}\nu\|_{L^\infty(\R^d)}\leq1$ 
for $0\leq j\leq d-1$ and $|\nu(B(x,r))|\leq r^s$
for all $x\in\R^d$, $r>0$ (in case $\nu$ is a distribution the latter condition 
should be reformulated appropriately). Is $\wt\gamma_s(E)\approx\gamma_s(E)$? 
It is shown in \cite{MPV-scalar} that the answer is affirmative for $s=1$ and negative for $0<s<1$, while it is unknown when $s>1$.

\bibliographystyle{alpha}
\bibliography{./refer}

\end{document}